\documentclass[reqno,fleqn,11pt]{amsart}

\usepackage{amsmath}
\usepackage{amssymb}
\usepackage{amsfonts}
\usepackage{amsthm}

\newcommand{\C}{\mathbb C}
\newcommand{\N}{\mathbb{N}}

\newcommand{\Z}{\mathbb{Z}}
\newcommand{\E}{\mathbb{E}}
\newcommand{\f}{\rightarrow}

\newcommand{\inj}{\operatorname{inj}}
\newcommand{\sys}{\operatorname{sys}}
\newcommand{\sgl}{\operatorname{sgl}}
\newcommand{\dias}{\operatorname{dias}}
\newcommand{\diam}{\operatorname{Diam}}
\newcommand{\Vol}{\operatorname{Vol}}
\newcommand{\Ent}{\operatorname{Ent}}
\newcommand{\Entalg}{\operatorname{Ent_{alg}}}
\newcommand{\Entvol}{\operatorname{Ent_{vol}}}

\newcommand{\GCD}{\operatorname{GCD}}
\newcommand{\id}{\operatorname{id}}

\newtheorem{thm}{Theorem}[section]
\theoremstyle{plain}

\newtheorem{lem}[thm]{Lemma}

\newtheorem*{ML}{Margulis Lemma}
\newtheorem*{LA}{Lemma A}
\newtheorem*{LB}{Lemma B}
\newtheorem*{LC}{Lemma C}
\newtheorem*{CBCG}{Corollaire 0.5 in \cite{BCG}}

\newtheorem*{comp}{Precompactness Theorem}
\newtheorem*{volume}{Volume estimate}
\theoremstyle{definition}
\newtheorem{defn}[thm]{Definition}

\theoremstyle{remark}
\newtheorem{rmk}[thm]{Remark}

\newtheorem*{case}{Case}
\newtheorem*{notation}{Notation}
\newtheorem{exmp}[thm]{Example}

\begin{document}

\title{Margulis Lemma, entropy and free products}
\author[F.Cerocchi]{Filippo Cerocchi}
\address{Dipartimento di Matematica "G.Castelnuovo", Universit\`{a} di Roma "Sapienza",
Piazzale Aldo Moro 5, 00185 Roma, Italy}
\address{Institut Fourier, Universit\'e Grenoble 1, 100 rue des maths, BP 74, 38402 St. Martin d'H\`eres, France}
\email{cerocchi@mat.uniroma1.it;   filippo.cerocchi@ujf-grenoble.fr; }

\begin{abstract}
We prove a Margulis' Lemma \textit{\`a  la} Besson Courtois Gallot,  for manifolds  whose fundamental group is a nontrivial free product $A*B$, without 2-torsion. Moreover, if $A*B$ is torsion-free we give a lower bound for the homotopy systole in terms of upper bounds on the diameter and the volume-entropy. We also provide examples and counterexamples showing the optimality of our assumption. Finally we give two applications of this result: a finiteness theorem and a volume estimate for reducible manifolds.
\end{abstract}

\maketitle

\section{Introduction}

The celebrated \textit{Margulis Lemma}, can be stated as follows:

\begin{ML}
Let $X$ be any compact Riemannian manifold of dimension $n\geq 2$, whose sectional curvature $\sigma(X)$ satisfies $-K^2\leq\sigma(X)<0$. Then:
$$\sup_{x\in X}\inj_x(X)\geq\frac{C_2(n)}{K},\quad\Vol(X)\geq\frac{C_1(n)}{K^n}$$
where $\inj_x(X)$ denotes the injectivity radius at $x$, and $C_1$, $C_2$ are two universal constants depending only on the dimension $n$.
\end{ML}

G. Besson, G. Courtois and S. Gallot gave in \cite{BCG} a more general version of Margulis' result: they replaced the strong assumption on the sectional curvature, by an algebraic hypothesis on the fundamental group ($\delta$\textit{-non ab\'elianit\'e}) together with an upper bound of the volume-entropy to obtain a lower bound for $l_x(X)$, the length of the shortest non-nullhomotopic geodesic loop based at some point $x$. We call this invariant the \textit{diastole} of $X$, $\dias(X)=\sup_{x}l_x(X)$ for easier reference throughout the paper:

\begin{CBCG}\label{becoga}
Let $\delta, H>0$. If $X$ is any Riemannian manifold whose fundamental group $\Gamma$ is $\delta$-nonabelian, and such that the commutation relation is transitive on $\Gamma\setminus\{\id\}$, with $\Entvol(X)\leq H$ we have:
$$\dias(X)\geq\frac{\delta\log(2)}{4+\delta}\cdot\frac{1}{H}.$$
\end{CBCG}
\medskip
Replacing $\sup_{x\in X}\inj_x(X)$ with $\dias(X)$ is the price to pay for dropping the negative curvature assumption.\medskip

We will denote by $\sys\pi_1(X)$ the homotopy systole of a compact Riemannian manifold $X$, \textit{i.e.} the length of the shortest non-contractible loop in $X$. We remark that if $\tilde X$ is the Riemannian universal covering of $X$ and $\tilde d$ its distance function we have
$\sys\pi_1(X)=\inf_{\tilde x\in\tilde X}\inf_{\gamma\in\pi_1(X)\setminus\{\id\}}\tilde d(\tilde x,\gamma\,\tilde x)$ and $\dias(X)=\sup_{\tilde x\in \tilde X}\inf_{\gamma\in\pi_1(X)\setminus\{\id\}}\tilde d(\tilde x,\gamma\,\tilde x)$\medskip\\
Following \cite{Wall} we say that a fundamental group is decomposable if it is isomorphic to a non trivial free product. We will say that a discrete group $\Gamma$ is without $2$-torsion (or $2$-torsionless) if there is no element $\gamma\in\Gamma$ such that $\gamma^2=\id$.\medskip\\
The main results in this paper are the following:

\begin{thm}\label{MT1}
Let $H>0$ and let $X$ be a connected Riemannian $n$-manifold such that $\Ent(X)\leq H$, whose fundamental group is decomposable, without 2-torsion. Then:
$$\dias(X)\geq\frac{\log(3)}{6H}.$$
\end{thm}

\begin{thm}\label{MT2}
Let $H,D>0$.
Let  $X$ be a compact Riemannian manifold such that $\Ent(X)\leq H$, $\diam(X)\leq D$, whose fundamental group is decomposable and torsion-free. Then we have:
$$\sys\pi_1(X)\geq\frac{1}{H}\cdot\log\left(1+\frac{4}{e^{2DH}-1}\right).$$
\end{thm}

The first theorem is based on the \textit{Kurosh subgroup theorem} (see \cite{Wall}, Theorem 3.1, pg. 151), the investigation of subgroups generated by "small elements" of $\pi_1(X)$ and a connectedness argument. The second one is a consequence of the \textit{Kurosh subgroup theorem} and of an estimate of the entropy of the Cayley graph of a free group generated by two elements when the corresponding edges have two different lengths. We remark that in \cite{BCG} there is an analogous statement which is valid only for $\delta$-thick groups (see \S5).\\
As a byproduct of Theorem \ref{MT2}, using a theorem of S. Sabourau (\cite{Sab2}, Proposition 3.5.4), under some extra geometric assumption, we obtain a precompactness and homotopy finiteness theorem:

\begin{comp}
Let $\mathfrak M_n^{dec}(D,V,H; l)$ denote the family of compact, Riemannian $n$-manifolds whose fundamental groups are decomposable and torsion-free, whose diameter, volume and volume-entropy are smaller than $D, V, H$ respectively, and such that the length of the shortest geodesic loop in the universal covering is greater than $l$. This family is precompact with respect to the Gromov-Hausdorff topology. Moreover there are at most finitely many distinct homotopy types in $\mathfrak M_n^{dec}(D,V,H;l)$.
\end{comp}

 Moreover, combining Theorem \ref{MT2} and the celebrated Isosystolic Inequality of Michael Gromov (\cite{Gromov3}, Theorem 0.1.A) we will prove a volume estimate, without curvature assumptions, for a certain class of Riemannian manifolds:

\begin{volume}
For any connected and compact, $1$-essential Riemannian $n$-manifold, $X$ with decomposable, torsion-free fundamental group and whose volume-entropy and diameter are bounded above by $H, D>0$ respectively, we have the following estimate:
$$\Vol(X)\geq\frac{C_n}{H^n}\cdot\log\left(1+\frac{4}{e^{2DH}-1}\right)^n$$
where $C_n>0$ is a universal constant depending only on the dimension $n$ (an explicit -although not optimal- upper bound to $C_n$ can be found in \cite{Gromov3}, Theorem 0.1.A).
\end{volume}

In section 2 we recall some basic facts about entropy. Section 3 is devoted to the proof of Theorem \ref{MT1}, while in section 4 we give the proof of Theorem \ref{MT2}, with some applications.
In section 5 we give examples showing that the class of manifolds covered by Theorem \ref{MT1} is distinct from the class considered in \cite{BCG}, Corollaire 0.5 (and in particular is orthogonal to the class of $\delta$-thick groups, cf. section 5). We also produce counterexamples showing that the torsion-free assumption in Theorem \ref{MT2} cannot be dropped: namely, we construct a manifold  $X$ with $\pi_1(X)=\Z_p*G$, for a non trivial group $G$, and a sequence of metrics with diameter and volume entropy bounded from above, whose homotopy systole tends to zero.

\section{Notations and background}

\begin{defn}\label{EDMG}
Let $(\Gamma,d)$ be a metric discrete group, \textit{i.e.} a group $\Gamma$ endowed with a left invariant distance such that $\#\{\gamma\;|\; d(\gamma,\id)<R\}<+\infty$, $\forall R>0$ (we call such a distance an admissible distance).
We define the entropy of the metric discrete group $(\Gamma,d)$, $\Ent(\Gamma,d)$ as:
$$\Ent(\Gamma,d)=\lim_{R\f\infty}\frac{1}{R}\log(\#\{\delta\;|\; d(\gamma,\delta)<R\})$$
This limit exists and does not depend on the element $\gamma$. 
\end{defn}

\begin{rmk}\label{ADGD}
We are interested in two different kinds of admissible distances on $\Gamma$. If $\Gamma$ is a finitely generated group and $\Sigma$ is a finite generating set we denote by $d_\Sigma$ the algebraic distance on $\Gamma$ associated to $\Sigma$. If $\Gamma$ is the fundamental group of a Riemannian manifold $X$, for any point $\tilde x$ in the Riemannian universal covering $\tilde X$, we define the admissible distance $d_{geo}$ on $\Gamma$ by: $d_{geo}(\gamma,\delta)=\tilde d(\gamma\,\tilde x,\delta\,\tilde x)$ where $\tilde d$ is the Riemannian distance on $\tilde X$.
\end{rmk}

\begin{defn}\label{algent}
Let $\Gamma$ be a discrete, finitely generated group. The algebraic entropy of $\Gamma$ is 
$\Entalg(\Gamma)=\inf_{\Sigma}\Ent(\Gamma,d_\Sigma)$,
where the $\inf$ is taken over the finite generating sets of $\Gamma$.
\end{defn}

\begin{defn}\label{EM}
Let $X$ be any Riemannian manifold. Its entropy is defined as:
$$\Ent(X)=\lim_{R\f\infty}\frac{1}{R}\log(\Vol(B(\tilde x,R)))$$
where $B(\tilde x,R)$ denotes the geodesic ball of radius $R$ in the Riemannian universal covering $\tilde X$ centered at $\tilde x$. The limit exists and it is easy to see that this does not depend on the point $\tilde x$. We remark that in the compact case this is just the volume-entropy of the Riemannian manifold.
\end{defn}

\begin{notation}
When we need to stress the dependence of $\Ent(X)$ from the Riemannian metric $g$ on $X$ we use the notation $\Ent(X,g)$ (or $\Entvol(X,g)$, in the compact case).
\end{notation}

We shall use the following basic properties of the entropy:
\begin{itemize}
\item[(0)] When $X$ is a Riemannian manifold and $\Gamma$ its fundamental group, then:
$\Ent(X)\geq\Ent(\Gamma, d_{geo})$ (\cite{BCG}, Lemma 2.3). Equality holds when $X$ is compact (see \cite{Rob2}, Proposition 1.4.7).
\item[(1)] Let $d_1\leq d_2$ be two admissible distances on $\Gamma$, then we have:\newline $\Ent(\Gamma,d_1)\geq\Ent(\Gamma,d_2).$
\item[(2)]  Let $d$ be an admissible distance on $\Gamma$ and let $\lambda>0$, then we have: $\Ent(\Gamma,\lambda d)=\frac{1}{\lambda}\Ent(\Gamma,d).$
\end{itemize}

\section{Proof of Theorem \ref{MT1}}

The proof of the Theorem is by contradiction and essentially relies on the following results:
\begin{itemize}
\item[i)] A structure theorem for finitely generated subgroups of free products (the well known \textit{Kurosh subgroup theorem}).
\item[ii)] The existence of a universal lower bound for the algebraic entropy of nontrivial free products $\neq\Z_2*\Z_2$ (\cite{delH}, \S VII.18):
\begin{equation}\label{AEE}
\Entalg(A*B)\geq\frac{\log(3)}{6}
\end{equation}
\item[iii)] The comparison between entropies of $X$ and $(\pi_1(X),d_{geo})$ (see property (0) of the entropy of a metric discrete group).
\end{itemize}\medskip
We recall that $\pi_1(X)$ is decomposable and $2$-torsionless. Let $l_0=\frac{\log(3)}{6H}$ and define the following family of sets:
$$\mathcal I(\tilde x,l_0)=\{\gamma\in A*B\setminus\{\id\}\;|\;\tilde d(\tilde x,\gamma\tilde x)<l_0\},\quad\forall\tilde x\in\tilde X.$$
Since $A*B$ acts by isometries on $\tilde X$ and the action is free and properly discontinuous, all these sets are finite. Moreover we underline the fact that they are symmetric (\textit{i.e.} if $\gamma\in\mathcal I(\tilde x,l_0)$, then $\gamma^{-1}\in\mathcal I(\tilde x, l_0)$).

\subsection{Three Lemmas} We will resume in the following three Lemmas the principal properties of the sets $\mathcal I(\tilde x,l_0)$:

\begin{LA}
For any $\tilde x\in\tilde X$, then $\mathcal I(\tilde x,l_0)$
\begin{itemize}
\item[(i)] either is included in $\gamma_{\tilde x} A\gamma_{\tilde x}^{-1}$, for at least one $\gamma_{\tilde x}\in A*B$;
\item[(ii)] or is included in $\gamma_{\tilde x} B\gamma_{\tilde x}^{-1}$, for at least one $\gamma_{\tilde x}\in A*B$;
\item[(iii)] or $\langle\mathcal I(\tilde x,l_0)\rangle\cong\Z$ and does not satisfy \textup{(i)} or \textup{(ii)}.
\end{itemize}
\end{LA}

\textbf{Proof.} Assume that conditions (i), (ii), (iii) are not verified; by the Kurosh Subgroup Theorem we know that the subgroup generated by $\mathcal I(\tilde x,l_0)$ writes
$$\langle\mathcal I(\tilde x,l_0(H))\rangle=C_1*\cdot\cdot *C_k*\gamma_1A_1\gamma_1^{-1}*\cdot\cdot *\gamma_r A_r\gamma_r^{-1}*\delta_1 B_1\delta_1^{-1}*\cdot\cdot *\delta_s B_s\delta_s^{-1}$$
where the $C_i$'s are infinite cyclic subgroups of $A*B$ which are not contained in any conjugate of $A$ or $B$, where $\gamma_i\neq\gamma_j$ and $\delta_i\neq\delta_j$ for $i\neq j$, and where $A_j$, $B_i$ are respectively subgroups of $A$, $B$. Since conditions (i), (ii), (iii) are not verified there should be at least two factors giving a nontrivial free product $\neq\Z_2*\Z_2$, hence by estimate \ref{AEE} we should have $\Entalg(\mathcal I(\tilde x,l_0)\rangle)\geq\log(3)/6$. On the other hand, by the triangle inequality, we have the inequality $d_{\mathcal I(\tilde x,l_0)}\cdot l_0> d_{geo}$, which is valid on $\langle\mathcal I(\tilde x,l_0)\rangle$; hence using properties (0), (1) and (2) of the entropy and the upper bound on the entropy of $X$ we prove that
$$\Ent(\langle\mathcal I(\tilde x, l_0)\rangle,d_{\mathcal I(\tilde x, l_0)})<\Ent\left(\langle\mathcal I(\tilde x,l_0)\rangle,\frac{1}{l_0}\cdot d_{geo}\right)=$$$$=l_0\cdot\Ent(\langle\mathcal I(\tilde x,l_0)\rangle, d_{geo})\leq H\cdot l_0=\frac{\log(3)}{6}$$
 which contradicts estimate (\ref{AEE}). $\Box$\newline

\begin{LB}
For all $\tilde x\in\tilde X$ there exists an $\varepsilon=\varepsilon(\tilde x)$ such that, for any $\tilde x'$, if $\tilde d(\tilde x',\tilde x)<\varepsilon$ the following inclusion holds:
$\mathcal I(\tilde x,l_0)\subseteq\mathcal I(\tilde x',l_0).$
\end{LB}

\textbf{Proof.} Let us fix
$\varepsilon<\frac{1}{2}\left[l_0-\sup_{\gamma\in\mathcal I(\tilde x,l_0)} \tilde d(\tilde x,\gamma\tilde x)\right]$
then by the triangular inequality we get the inclusion. $\Box$\newline

\begin{LC}
For all $\tilde x\in\tilde X$ and for all $\gamma\in\Gamma$ the following equality holds: $$\gamma\mathcal I(\tilde x,l_0)\gamma^{-1}=\mathcal I(\gamma(\tilde x), l_0).$$
\end{LC}

\textbf{Proof.}  Let $\delta\in\mathcal I(\tilde x,l_0)$; then $\gamma\delta\gamma^{-1}$ satisfies the inequality:
$$\tilde d(\gamma\tilde x,\gamma\delta\gamma^{-1}\cdot\gamma\tilde x)=\tilde d(\gamma\tilde x,\gamma\delta\tilde x)=\tilde d(\tilde x,\delta\tilde x)<l_0$$
hence $\gamma\delta\gamma^{-1}\in\mathcal I(\gamma\tilde x,l_0)$. To obtain the reverse inclusion suppose to have $\sigma\in\mathcal I(\gamma\tilde x, l_0)$, proceeding as before we obtain $\gamma^{-1}\sigma\gamma\in\mathcal I(\tilde x, l_0)$; hence $\sigma$ is a $\gamma$-conjugate  of an element in $\mathcal I(\tilde x,l_0)$. $\Box$\newline

\subsection{End of the proof}

\textit{We assume that $\mathcal I(\tilde x,l_0)\neq\varnothing$}, for all $\tilde x\in\tilde X$ and we shall show that this leads to a contradiction.
Let us now define the following sets:
\begin{itemize}
\item $\tilde X_1=\{\tilde x\in\tilde X\;|\;\exists\gamma\in A*B\mbox{ such that } \mathcal I(\tilde x, l_0)\subseteq \gamma A\gamma^{-1}\}$;
\item $\tilde X_2=\{\tilde x\in\tilde X\;|\;\exists\gamma\in A*B\mbox{ such that } \mathcal I(\tilde x, l_0)\subseteq \gamma B\gamma^{-1}\}$;
\item $\tilde X_3=\{\tilde x\in\tilde X\setminus(\tilde X_1\cup\tilde X_2)\;|\;\exists\tau\in A*B\mbox{ such that } \langle\mathcal I(\tilde x,l_0)\rangle=\langle\tau\rangle\}$;
\end{itemize}
the next lemma enlightens some key properties of these sets:

\begin{lem}\label{opendisjoint}
The $\tilde X_i$'s are open and disjoint.
\end{lem}

\textbf{Proof.} $\tilde X_3$ is disjoint from $\tilde X_1$ and $\tilde X_2$ by definition, whereas $\tilde X_1$ and $\tilde X_2$ are disjoint since $\mathcal I(\tilde x, l_0)\neq\varnothing$, $\id\not\in\mathcal I(\tilde x,l_0)$ and $$\gamma A\gamma^{-1}\cap\delta B\delta^{-1}=\{\id\}\;\; \forall\gamma,\delta\in A*B.$$
Now we will prove that $\tilde X_i$ is open. Let us take a point $\tilde x$ in $\tilde X_i$; by Lemma B, for $\tilde x'$ in an open neighbourhood of $\tilde x$ we have the inclusion, $\mathcal I(\tilde x,l_0)\subseteq\mathcal I(\tilde x',l_0)$. Now, by Lemma A, for $\tilde x'$ one condition between (i), (ii) and (iii) should hold, \textit{i.e.} $\exists j\in\{1,\,2,\, 3\}$ such that $\tilde x'\in\tilde X_j$. As the $\tilde X_i$ are disjoint, by the inclusion above it follows that if $\tilde x\in\tilde X_i$, then also $\tilde x'\in\tilde X_i$. Hence the $\tilde X_i$'s are open. $\Box$\newline

Since the $\tilde X_i$'s are open and disjoint subsets of $\tilde X$ and $\tilde X$ is connected one of the following conditions should hold:
\begin{enumerate}
\item $\tilde X=\tilde X_1$;
\item $\tilde X=\tilde X_2$;
\item $\tilde X=\tilde X_3$;
\end{enumerate}
we will now show that each of these conditions leads to a contradiction.

\begin{case}[1]
We shall prove that there exists $\gamma_0$, independent from $\tilde x$ such that all  the sets $\mathcal I(\tilde x,l_0)$ belong to the same conjugate $\gamma_0A\gamma_0^{-1}$ of $A$ in $A*B$. For each fixed $\hat\gamma\in (A*B)/A$ we define the subset of $\tilde X_1$:
$$\tilde X_1(\hat\gamma)=\{\tilde x,\in\tilde X\;|\;\exists\gamma\in\hat\gamma\mbox{ such that }\mathcal I(\tilde x,l_0)\subseteq \gamma A\gamma^{-1}\}$$
and we remark that since $\tilde X=\tilde X_1$, we have $\tilde X=\cup_{\hat\gamma}\tilde X_1(\hat\gamma)$.  The sets $\tilde X_1(\hat\gamma)$ are disjoint: the proof is analogous to the proof of Lemma \ref{opendisjoint}. Moreover every 
$\tilde X_1(\hat\gamma)$ is open: let $\tilde x\in\tilde X_1(\hat\gamma)$ and consider a $\tilde x'$ at distance $\tilde d(\tilde x,\tilde x')<\varepsilon$, where $\varepsilon$ is chosen as in Lemma B. By Lemma B we know that $\mathcal I(\tilde x,l_0)\subseteq\mathcal I(\tilde x',l_0)$ so if $\tilde x\in\tilde X_1(\hat\gamma)$ since by assumption $\mathcal I(\tilde x,l_0)\neq\varnothing$ then also $\tilde x'\in\tilde X_1(\hat\gamma)$, because $\tilde X_1(\hat\gamma)$ and $\tilde X_1(\hat\gamma')$ are disjoint if $\hat\gamma\neq\hat\gamma'$. Hence $\tilde X$ is covered by the family of disjoint, open sets $\tilde X_1(\hat\gamma)_{\hat\gamma\in(A*B)/A}$, and by connectedness of $\tilde X$, there exists one $\hat\gamma_0\in(A*B)/A$ such that $\tilde X=\tilde X_1(\hat\gamma_0)$. So there exists $\gamma_0\in A*B$ such that every subset $\mathcal I(\tilde x,l_0)$ is included in some $\gamma A\gamma^{-1}$ for $\gamma\in\hat\gamma_0=\gamma_0A$ (\textit{i.e.} in $\gamma_0 A\gamma_0^{-1}$). It follows that the whole subgroup $G=\langle\mathcal I(\tilde x,l_0)\rangle_{\tilde x\in\tilde X}$ is included in $\gamma_0 A\gamma_0^{-1}$. By construction of $G$ and by Lemma C, $G$ should be a normal subgroup of $A*B$, and this is a contradiction, as no normal subgroup of a nontrivial free product is included in a conjugate of one factor. 
\end{case}

\begin{case}[2]
 Proof of Case (2) is analogous to Case (1).  
\end{case}

\begin{case}[3]
Let $T\subset A*B$ be the subset of primitive elements of $A*B$ \footnote{\textit{i.e.} elements which cannot be written as powers of any other element in $A*B$.} with infinite order, which are not contained in any conjugate of $A$ or $B$. Similarly to (1) we will show that there exists $\tau_0\in T$ such that $\langle\mathcal I(\tilde x, l_0)\rangle\subset\langle\tau_0\rangle$ for all $\tilde x\in\tilde X$. For each $\tau\in T$ let
$$\tilde X_3(\tau)=\tilde X_3(\tau^{-1})=\{\tilde x\;|\;\exists k\in\Z\mbox{ such that }\langle\mathcal I(\tilde x,l_0)\rangle=\langle\tau^k\rangle\}$$
We want to show that  $\tilde X=\bigcup_{\tau\in T}\tilde X_3(\tau)$; by assumption every subgroup $\langle\mathcal I(\tilde x,l_0)\rangle$ is isomorphic to an infinite cyclic subrgoup $\langle\gamma\rangle$ hence it suffices to show that for any element $\gamma$ in $A*B$ there exists an element $\tau\in T$ and $k\in\Z$ such that $\gamma=\tau^k$. We argue by contradiction: assume that there is an element $\gamma$, which cannot be written as a power of a primitive element $\tau\in T$; by definition there exists a sequence of elements $\{\gamma_n\}_{n\in\N}$ and a sequence of integers $\{p_n\}_{n\in\N}$ such that $\gamma_0=\gamma$ and $\gamma_i=(\gamma_{i+1})^{p_{i+1}}$ (with $|\prod_1^{i}p_j|\f\infty$). Let $\Sigma$ be any generating system and $d_0=d_\Sigma(\gamma,\id)$; for any $i\in\N$ consider $\delta_i$ a cyclically reduced word associated to $\gamma_i$ and let $N_1(i)=\frac{l_\Sigma(\gamma_i)-l_\Sigma(\delta_i)}{2},\, N_2(i)=l_\Sigma(\delta_i)$. Observe that since $\gamma\neq\id$ we shall have $N_2(i)\geq 1$. Then $d_{\Sigma}(\gamma_i^{p_i},\,\id)=2N_1(i)+|p_i|\cdot N_2(i)$ so that:
$$d_0=d_\Sigma(\gamma_{i}^{\prod_1^{i}p_j},\id)\geq|\prod_1^i p_j|\cdot N_2(i)\geq |\prod_1^ip_j|$$
which gives a contradiction for $i\f\infty$.\smallskip\\
Let us show that the sets $\tilde X_3(\tau)$ are disjoint. Actually assume that $\tilde x\in \tilde X_3(\tau)\cap\tilde X_3(\tau')$, for  $\tau'\neq\tau^{\pm1}$; then there exists $k,\,k'$ such that $\sigma=\tau^k=(\tau')^{k'}$ generates $\langle\mathcal I(\tilde x, l_0)\rangle$. The following lemma shows that one among $\tau$ and $\tau'$ is not primitive, a contradiction:

\begin{lem}[Primitive powers in free products]\label{PPFP}
Let $\gamma,\,\gamma'\in A*B\setminus\{\id\}$ be such that $\gamma^s=(\gamma')^{s'}$ ($s,\, s'\in\N$), and assume that they are not contained in any conjugate of $A$ or $B$. Then there exists an element $\tau\in A*B$ and $q,\,q'\in\N$ such that $\gamma=\tau^q,\, \gamma'=\tau^{q'}$. In particular if $\gamma$, $\gamma'$ are primitive elements then $q,\,q'\in\{-1,\,1\}$ and $\gamma=(\gamma')^{\pm1}$.
\end{lem}

The proof of this  Lemma is rather simple but tedious and will be given in the Appendix.  We now prove that the sets $\tilde X_3(\tau)$ are open. Let $\tilde x\in\tilde X_3(\tau)$: we know by Lemma B that, for all $\tilde x'$ sufficiently close to $\tilde x$ we have $\mathcal I(\tilde x,l_0)\subseteq\mathcal I(\tilde x',l_0)$; if $\tilde x'\not\in\tilde X_3(\tau)$ then
$$\langle(\tau')^{k'}\rangle=\langle\mathcal I(\tilde x',l_0)\rangle\supset\langle\mathcal I(\tilde x,l_0)\rangle=\langle\tau^k\rangle$$
for some $k,\, k'\in\Z$ and $\tau'\in T$ different from $\tau^{\pm 1}$. Again Lemma \ref{PPFP} implies that $\tau'$ or $\tau$ is not primitive, a contradiction.
It follows, by connectedness, that  $\tilde X=\tilde X_3(\tau)$ for some  fixed $\tau\in T$. Therefore, the group $\langle\tau\rangle$ (that contains $\langle\mathcal I(\tilde x,l_0)\rangle$, for any $\tilde x\in\tilde X$) is a normal subgroup: in fact, for any $\gamma\in A*B$ there exists $k,\, k'$ such that $\langle\mathcal I (\gamma\,\tilde x,l_0)\rangle=\langle\tau^k\rangle=\langle\gamma\tau^{k'}\gamma^{-1}\rangle$; since $\tau$ and $\gamma\tau\gamma^{-1}$ are both primitive elements it follows that $\gamma\tau\gamma^{-1}=\tau^{\pm 1}$. Hence $\langle\tau\rangle$ is an infinite cyclic subgroup of $A*B$, which is normal. This is not possible since no free product different from $\Z_2*\Z_2$ admits an infinite cyclic normal subgroup\footnote{In fact by \cite{Wall}, Theorem 3.11, p. 160, every normal subgroup in $A*B$ must have finite index, and an infinite cyclic group in $A*B$ can have finite index  if and only if $A*B=\Z_2*\Z_2$.}.
This excludes also Case (3).
\end{case}
Therefore $\mathcal I(\tilde x,l_0)=\varnothing$ for some $\tilde x\in\tilde X$, which proves Theorem \ref{MT1}. $\Box$

\section{Proof of Theorem \ref{MT2} and Applications}

\subsection{Proof of Theorem \ref{MT2}} Let $X$ be a compact Riemannian manifold with decomposable, torsion free fundamental group and assume the bounds $\Ent(X)\leq H$, $\diam(X)\leq D$. The proof relies on the following Lemma:

\begin{lem}\label{entsemigroup}
Let $\mathcal C(\Gamma,\{\gamma_1,\,\gamma_2\})$ be the Cayley graph of a free group with two generators $\gamma_1,\,\gamma_2$. Let $d_l$ be the left invariant distance on $\mathcal C(\Gamma,\{\gamma_1,\,\gamma_2\})$, defined by the conditions $d_l(\id,\gamma_1)=l(\gamma_1)$ and $d_l(\id,\gamma_2)=l(\gamma_2)$. Then $h=\Ent(\Gamma,d_l)$ solves the equation: 
\begin{equation}\label{entest}
(e^{h\cdot l(\gamma_1)}-1)(e^{h\cdot l(\gamma_2)}-1)=4
\end{equation}
\end{lem}

\textbf{Proof.} Let $I(c)=\sum_{\gamma\in\Gamma} e^{-c\, d_l(\id,\gamma)}$. Since $l(\gamma_1)$ and $l(\gamma_2)$ are strictly positive, the entropy of $(\Gamma,d_l)$ is finite (but not necessarily bounded independently from $l(a)$, $l(b)$), $\{c>0\;|\; I(c)<+\infty\}\neq\varnothing$ and $\Ent(\Gamma,d_l)=\inf\{c>0\;|\; I(c)<+\infty\}$. Let us define the sets $S_{\gamma_1^{\pm 1}},\, S_{\gamma_2^{\pm 1}}$ as the sets of elements of $\Gamma$ whose reduced writing starts by $\gamma_1^{\pm 1}$, $\gamma_2^{\pm 1}$ (respectively). We define $I_s(c)=\sum_{\gamma\in S_s} e^{-c\,d_l(\id,\gamma)}$
where $s\in\{\gamma_1,\,\gamma_1^{-1},\,\gamma_2,\,\gamma_2^{-1}\}$. By definition we have $I(c)=1+I_{\gamma_1}(c)+I_{\gamma_1^{-1}}(c)+I_{\gamma_2}(c)+I_{\gamma_2^{-1}}(c)$ that is:
\begin{equation}\label{I}
I(c)=1+2(I_{\gamma_1}(c)+I_{\gamma_2}(c))
\end{equation}
Moreover, since $S_{\gamma_1}=\gamma_1\cdot(S_{\gamma_1}\cup S_{\gamma_2}\cup S_{\gamma_2^{-1}})$ we have: $I_{\gamma_1}(c)= e^{-c\, l(\gamma_1)}\cdot(I_{\gamma_1}(c)+I_{\gamma_2}(c)+I_{\gamma_2^{-1}}(c))$. Hence we have: $I_{\gamma_1}(c)+e^{-c\,l(\gamma_1)}\,I_{\gamma_1^{-1}}(c)=e^{-c\, l(\gamma_1)}\cdot I(c)$ and since $I_{\gamma_1}(c)=I_{\gamma_1^{-1}}(c)$ we get:
\begin{equation}\label{Igamma1}
I_{\gamma_1}(c)=\frac{I(c)}{(e^{c\, l(\gamma_1)}+1)}
\end{equation}
Analogously one has:
\begin{equation}\label{Igamma2}
I_{\gamma_2}(c)=\frac{I(c)}{(e^{c\,l(\gamma_2)}+1)}
\end{equation}
Now we plug equations (\ref{Igamma1}) and (\ref{Igamma2}) into equation (\ref{I}):
$$I(c)=1+2\cdot\left[\frac{1}{e^{c\,l(\gamma_1)}+1}+\frac{1}{e^{c\,l(\gamma_2)}+1}\right]\cdot I(c)$$
and since $I(c)\f+\infty$ as $c\f h_+$ we see that equation (\ref{entest}) holds. $\Box$ \medskip\\
\textbf{End of the proof of Theorem \ref{MT2}.} Let us now fix a point $\tilde x\in\tilde X$; let $\sigma_1,\,\sigma_2\in A*B$ be two elements such that $\langle\sigma_1,\sigma_2\rangle\simeq\mathbb F_2$ is a nontrivial free product (hence a free group since $A*B$ is torsion free).  Let us denote $l(\sigma_1)=\tilde d(\sigma_1\tilde x,\tilde x)$, $l(\sigma_2)=\tilde d(\sigma_2\tilde x,\tilde x)$. Using Lemma \ref{entsemigroup}, we obtain

$$H\geq\Ent(\pi_1(X), d_{geo})\geq\Ent(\langle\sigma_1,\sigma_2\rangle, d_{geo})\geq\Ent(\langle\sigma_1,\sigma_2\rangle, d_l)\geq$$
$$\geq\frac{1}{l(\sigma_1)}\cdot\log\left(1+\frac{4}{e^{\Ent(\langle\sigma_1,\sigma_2\rangle, d_l)\,l(\sigma_2)}}\right)\geq\frac{1}{l(\sigma_1)}\cdot\log\left(1+\frac{4}{e^{H\, l(\sigma_2)}}\right)$$
from which we deduce:
\begin{equation}\label{LSE}
\tilde d(\tilde x,\sigma_1\tilde x)\geq\frac{1}{H}\cdot\log\left(1+\frac{4}{e^{H\,\tilde d(\tilde x,\sigma_2\tilde x)}}\right)
\end{equation}
Let $\sigma$ be a geodesic loop realizing $\sys\pi_1(X)$ and let $\tilde x$ belong to $\sigma$. Let $\Sigma=\{\tau_i\}$ be a finite generating set such that $\tilde d(\tau_i\tilde x,\tilde x)\leq 2D$ (\cite{Gromov}, Proposition 5.28). There exists at least one $\tau_i\in\Sigma$ such that $\langle\tau_i,\sigma\rangle\simeq\mathbb F_2$ is a free product (hence a free group), since $\Sigma$ is a generating set and $A*B$ is a free product. As $\tilde d(\tilde x,\tau_i\tilde x)\leq 2D$, the inequality (\ref{LSE}) applied to $\sigma$, $\tau_i$ gives:
$$\sys\pi_1(X)=\tilde d(\sigma\tilde x,\tilde x)\geq\frac{1}{H}\cdot\log\left(1+\frac{4}{e^{2DH}-1}\right)$$
This ends the proof of Theorem \ref{MT2}. $\Box$

\begin{rmk}
Other lower bounds for the homotopy systole, with upper bounds for the volume entropy and the diameter (in addition to some algebraic assumption on $\pi_1(X)$) have been proved in \cite{BCG}.
However in the next section we shall show a quite large class of examples where our estimate can be applied but not those of \cite{BCG}.
\end{rmk}

\subsection{Applications}
Let $Y$ be a complete Riemannian manifold; we will denote by $\sgl(Y)$the length of the shortest  (possibly homotopically trivial) geodesic loop in $Y$.

\begin{thm}\label{quotients}
Let $X$ be a simply connected, Riemannian manifold. The family $\mathfrak M_{X}^{dec}(D,V,H)$ of compact, Riemannian quotients of $X$ with torsionless, decomposable fundamental group such that diameter, volume and volume-entropy are bounded by $D, V, H$, respectively, is finite up to homotopy.\end{thm}

\textbf{Proof.} It is a direct consequence of the Precompactness Theorem.$\Box$\newline

\textbf{Proof of the Precompactness Theorem.} We follow the proof of the Proposition 4.3 of \cite{BCG}. Let $X\in\mathfrak M_n^{dec}(D,V,H;l)$; and let $\tilde X$ be its Riemannian universal covering. For any $x\in X$, the distance between distinct points $\tilde x_1,\tilde x_2$ of the $x$-fiber in $\tilde X$ is greater or equal to $\sys\pi_1(X)$, so $B(x,\frac{\sys\pi_1(X)}{2})$ is isometric to $\tilde B(\tilde x,\frac{\sys\pi_1(X)}{2})$, for $\tilde x$ in the $x$-fiber in $\tilde X$. Since in $B(x,\frac{\sys\pi_1(X)}{2})$ we do not have geodesic loops of length less than $\sgl(\tilde X)$  (by definition of $\sgl(\tilde X)$), it follows that $\sgl(X)=\min\{\sys\pi_1(X),\sgl(\tilde X)\}$. Hence by Theorem \ref{MT2} and by the assumption we made on $\sgl(\tilde X)$ it follows that $\sgl(X)\geq l_0=\min\{ l,\frac{1}{H}\cdot\log(1+\frac{4}{e^{2DH}-1})\}$. A theorem of Sabourau (\cite{Sab}, Theorem A) states that if $M$ is a complete Riemannian manifold of dimension $n$ there exists a constant $C_n$, depending only on the dimension of $M$ such that $\Vol(B(x,R))\geq C_nR^n$, for every ball of radius $R\leq\frac{1}{2}\sgl(M)$. This means that we can bound the maximum number $N(X,\varepsilon)$ of disjoint geodesic balls in $X$ of radius $\varepsilon$ by the function $V/(C_n\varepsilon^n)$, and the estimate holds for any manifold in $\mathfrak M_n^{dec} (D,V,H;l)$ (obviously for $\varepsilon\leq l_0/2$). Then using the Gromov's packing argument as shown in \cite{Gromov} \S 5.1-5.3 and, for example, in \cite{Fuk}, Lemma 2.4, we get the precompactness of the family $\mathfrak M_n^{dec}(D,V,H;l)$. For what concerns finiteness of homotopy types in $\mathfrak M_n^{dec}(D,V,H;l)$ it follows from the precompactness of the family and from another result by Sabourau (\cite{Sab2}, Proposition 3.5.4), which states that if $M, N$ are two Riemannian $n$-manifolds satisfying $\sgl\geq l_0$ and such that their Gromov-Hausdorff distance $d_{GH}(M, N)$ is less than $\beta_n l_0$, then they have the same homotopy type. $\Box$\\

\begin{rmk}\label{comparison}
We want to compare our precompactness and finiteness theorem with the classical ones by J. Cheeger and M. Gromov (\cite{Ch}, \cite{Gromov3} \S 8.20, \cite{Gromov4}) and with more recent results by I. Belegradek (\cite{bel}). The first finiteness result has been obtain combining the results of Cheeger and Gromov (see \cite{Fuk}, Theorem 14.1): they considered the class of $n$-manifolds with bounded sectional curvature $|\sigma|\leq 1$, volume bounded below by a universal constant $v>0$ and diameter bounded above by a constant $D>0$ and they proved the finiteness of diffeomorphism classes (the proof given by Gromov uses the Lipschitz precompactness of the family and his rigidity theorem, see also \cite{Kat}, \cite{Fuk}). Observe that the assumptions of this result implies our geometric assumptions: in fact the first assumption (on curvature) implies the boundedness of the volume-entropy, while the three assumptions both imply the boundedness of the volume and a lower bound of $\sgl$. On the other hand in our case we can only achieve the finiteness of homotopy types.

Another finiteness theorem has been proved by Gromov in \cite{Gromov4}: he assumes to have the following bounds $-1\leq\sigma<0$ on the sectional curvature and an upper bound for the volume, $V$, and establishes the finiteness of diffeomorphism types for the class of Riemannian manifolds (of dimension $n\neq 3$) satisfying these bounds. The result is a consequence of Theorem 1.2 in \cite{Gromov4}, which gives an upper bound for the diameter of a negatively curved Riemannian manifold of sectional curvature $-1\leq\sigma<0$ in terms of its volume, combined with Cheeger's finiteness theorem and with Margulis' Lemma. We observe that also in this case since the Riemannian universal covering $\tilde X$ satisfies $\sgl(\tilde X)=+\infty$, the assumptions made by M. Gromov imply our geometric assumptions; moreover, the prescribed sign and the boundedness of the sectional curvature impose algebraic restrictions on the possible fundamental groups.

More recently I. Belegradek showed that once we fix a group $\Gamma$ for any $b\in[-1,0)$ there exist at most finitely many nondiffeomorphic closed Riemannian manifolds satisfying $-1\leq\sigma\leq b<0$ and whose fundamental group is isomorphic to $\Gamma$ (\cite{bel}, Corollary 1.4). Here the isomorphism class of the fundamental group is prescribed, but no assumption has been made on the volume and the diameter (whereas boundedness of sectional curvature implies the boundedness of the volume-entropy).\medskip
\end{rmk}

Another application of Theorem \ref{MT2} is a volume estimate for $1$-essential, compact, Riemannian $n$-manifolds with decomposable torsion free fundamental groups. We recall that a manifold $X$ is said to be $1$-essential whenever it admits a map $f$ into a $K(\pi,1)$-space $K$, such that the induced homomorphism $H_n(X,\Z)\f H_n(K,\Z)$ does not vanish.\newline

\textbf{Proof of the Volume estimate.} Just combine the estimate for the homotopy systole in Theorem \ref{MT2} with the inequality $\sys\pi_1(X)^n C_n\leq\Vol(X)$ proved in Theorem 0.1.A in \cite{Gromov3}. $\Box$

\section{Examples and Counterexamples}
In \cite{Zuddas} the following class of groups is defined:
 a $N$-nonabelian group is a group $\Gamma$ without nontrivial normal, abelian subgroups, such that the commutation relation is transitive on $\Gamma\setminus\{\id\}$ and such that  $\forall\gamma_1,\gamma_2\in\Gamma$ that do not commute, there exist two elements in $B(\id,N)\subset(\langle\gamma_1, \gamma_2\rangle, d_{alg})$ (here $d_{alg}$ denotes the algebraic distance $d_{\{\gamma_1,\gamma_2\}}$), which generate a  free semi-group (we call this last property the FSG($N$) property). This notion is inspired by the one of $\delta$-nonabelian group, introduced in \cite{BCG}. We remark that in general a $\delta$-nonabelian group (in the sense of \cite{BCG}) is not $N$-nonabelian (in the sense of \cite{Zuddas}), however $\delta$-nonabelian groups whose commutation relation is transitive are always $[\frac{4}{\delta}]$-nonabelian. Simple examples of $N$-nonabelian groups are:
\begin{itemize}
\item $\delta$-thick groups in the terminology of \cite{BCG} (\textit{i.e.} fundamental groups of Riemannian manifolds with sectional curvature less or equal to $-1$ and injectivity radius greater than $\delta$) are $[\frac{4}{\delta}]$-nonabelian.
\item  Free products of $\delta$-thick groups, and free products of $\delta$-thick groups with abelian groups are $[\frac{4}{\delta}]$-nonabelian.
\item More generally free products of $N$-nonabelian groups and free products of $N$-nonabelian groups with abelian groups are $N$-nonabelian (this is an easy corollary of Proposition 1.3 in \cite{Zuddas}).
\item $\pi_1(X)*\pi_1(Y)$, the free product of the fundamental groups of two compact Riemannian manifolds $X$, $Y$ with sectional curvature less or equal to $-1$ is $N$-nonabelian for $N\geq4\cdot\max\{\frac{1}{\inj(X)},\frac{1}{\inj(Y)}\}$ (see \cite{Zuddas}, \S1.4).
\item $\Z_n*\Z$, for every odd integer $n$, is  $N$-nonabelian for $N=4$ (again in \cite{Zuddas}, \S 1.4). This example is important since it shows that there are $N$-nonabelian groups which are not $\delta$-nonabelian.
\end{itemize}
\medskip
Following Zuddas we remark that $\delta$-nonabelian groups in the sense of \cite{BCG} satisfy a strictly stronger condition than the FSG($N$) property.
We give some examples of manifolds which satisfy the assumptions of our Theorem \ref{MT2}, whose fundamental groups are not $N$-nonabelian or $\delta$-nonabelian:

\begin{exmp}[Connected sums with flat manifolds] Consider $Y=X\# Z$ the connected sum of a quotient $X$ of $\E^n$ by the action of a discrete, nonabelian, torsion free and cocompact subgroup of $\mathrm{Is}(\E^n)$, with a compact manifold $Z$ whose fundamental group is torsion free and non trivial. Let $A=\pi_1(X)$ and $B=\pi_1(Z)$ then, if $n\geq 3$, $\pi_1(X\# Z)=A*B$. Then, the group $A*B$ does not possess the FSG($N$)-property. As $A$ is nonabelian there exist two elements $a_1, a_2$ which do not commute; since $A$ is a Bieberbach group, it is a group of polynomial growth (it contains a $\Z^{n}$ with finite index): this means that $A$ does not contain any free semigroup. So $A$ has a couple of elements not commuting and such that $\nexists N\in\N$ for which we can find two elements in $B(\id, N)\subset(\langle a_1, a_2\rangle, d_{\{a_1,a_2\}})$ that generate a free semigroup. 
\end{exmp}

\begin{exmp}[Connected sums with infranilmanifolds] 
More generally the above arguments hold for the connected sums with infranilmanifolds (a infranilmanifold is the quotient of a simply connected nilpotent Lie group by a nonabelian, torsion free, quasi-crystallographic group, see \cite{dek}, section 2.2). For instance if $X_k$ is the quotient of the Heisenberg group by
$$\Gamma_k=\langle a, b, c |\mbox{ }[b,a]=c^k,\mbox{ }[c,a]=[c,b]=\id\rangle$$
(where $a, b ,c$ are the standard generators) and $Y_k=X_k\# M$, where $M$ has non trivial, torsion free fundamental group, we have an infinite number of distinct differentiable manifolds, all non $N$-nonabelian, to which our Theorem \ref{MT2} applies (for \textit{any} choice of a Riemannian metric on $Y_k$).
\end{exmp}

\begin{rmk}
It is well known that connected sums of $1$-essential $n$-manifolds with other $n$-manifolds are still $1$-essential, so the examples above also provide a class of manifolds for which our Volume estimate holds.
\end{rmk}

Let us do some comments about Theorem \ref{MT1}. First of all we considerably enlarge the class of manifolds for which the Margulis Lemma \textit{\`a  la} Besson Courtois Gallot holds: in fact the only free products considered in \cite{BCG} and \cite{Zuddas} were free products of $N$-nonabelian groups or free proudcts of $N$-nonabelian groups with certain abelian groups; on the contrary we consider free products without restrictions, except for the $2$-torsionless assumption. Finally a remark about the necessity of requiring $\Gamma$ without 2-torsion: it might be sufficient to ask $\Gamma\neq\Z_2*\Z_2$ -\textit{i.e.} to exclude the unique case of a free product for which estimate (\ref{AEE}) does not hold-; however our proof of  the Theorem \ref{MT1} is not sufficient to conclude even in the case when $\Gamma=A*B$ and $B$ (or $A$) admits $2$-torsion.

\begin{exmp}
 We now exhibit a family of manifolds which proves that the 'torsion free' assumption of Theorem \ref{MT2} cannot be dropped. Fix  a $p\in\N$ and let $X$ be the connected sum of a lens space $M_p=S^3/\Z_p$ with any non simply connected manifold $Y$. We can endow $X$ with a family of metrics $g_\varepsilon$, such that for all $\varepsilon\in(0,1]$:
\begin{enumerate}
\item $\diam(X, g_\varepsilon)\leq D$ for a suitable $D\in(0,+\infty)$;
\item $\sys\pi_1(X,g_\varepsilon)\leq 2\pi\varepsilon/p$.
\item $\Ent(X,g_\varepsilon)\leq H$ for a suitable $H\in(0,+\infty)$;
\end{enumerate}
\medskip
5.4.1. \textit{Construction of the metrics $g_\varepsilon$ on $X$.} First we recall that if $S^3$ is endowed with the canonical metric (which we will denote by $h_1$ in the sequel) then we have an isometric action of $S^1$ which can be described as follows:
$$(S^3,can)=\{(u,v)\in\C^2 | |u|^2+|v|^2=1\},$$
$$S^1\times S^3\f S^3,\quad (e^{i\theta}, (u,v))\f(e^{i\theta} u, e^{i\theta} v).$$
Then $M_p=S^3/\langle e^{\frac{2\pi}{p}}\rangle$, so $\pi_1(M_p)=\langle \sigma_p\rangle$ where $\sigma_p=e^{\frac{2\pi}{p}}$. Let $\gamma_p$ be the shortest non contractible loop of $M_p$ representing $g_p$ (corresponding to an arc $\tilde\gamma_p$ on a maximal circle $\tilde\gamma$ in $S^3$). Since the normal bundle of $\tilde\gamma$, $N_\delta(\tilde\gamma)$,  is topologically trivial (\textit{i.e.} $N_{\delta}(\tilde\gamma)\simeq\tilde\gamma\times D_\delta^2$ where $D_\delta^2$ is a euclidean disk of radius $\delta$), in a tubular neighbourhood  $N_\delta(\tilde\gamma)$ we modify the canonical metric of $S^3$ only in the direction tangent to the $S^1$-action by a smooth factor $\lambda_\varepsilon(r)$ (where $r$ is the distance from the maximal circle $\tilde\gamma$), where $\lambda_\varepsilon\leq 1$ everywhere, $\lambda_\varepsilon(r)\equiv 1$ outside $N_\frac{2\delta}{3}(\tilde\gamma)$ and $\lambda_\varepsilon\equiv\varepsilon^2$ on $N_{\frac{\delta}{3}}(\tilde\gamma)$. We remark that $\lambda_\varepsilon$ can be constructed in $N_{\frac{2\delta}{3}}(\tilde\gamma)\setminus N_{\frac{\delta}{3}}(\tilde\gamma)$ in order to keep the sectional curvatures bounded below by a negative constant $C$, independent from $\varepsilon$. We obtain a new metric $\tilde h_\varepsilon$.
The new metric is still invariant by the action of $\Z_p$ defined before, hence the action of $\Z_p$ on $(S^3,\tilde h_\varepsilon)$ is still isometric and induce a metric $h_\varepsilon$ on $S^3/\Z_p$. Now, for any fixed metric $k$ on $Y$ we can glue $(M_p,h_\varepsilon)$ to $(Y,k)$ by gluing $M_p\setminus B_m$ and $Y\setminus B_y$ on the boundaries of two small balls $B_m$ and $B_y $ of $(M_p, h_\varepsilon)$ and $(Y,k)$ (respectively) such that $B_m$ lies outside $N_\delta(\gamma)$. We will call $(X,g_\varepsilon)$ the manifolds obtained in this way. \newline\newline
5.4.2. \textit{Proof of }(1), (2), (3). As the metric $g_\varepsilon$ is equal to $g_1$ except for the tubular neighbourhood of $\tilde\gamma$ and as $\lambda_\varepsilon$ is bounded above by $1$ we remark that $0<d_0<\diam(X,g_\varepsilon)\leq\diam(X,g_1)=D$ for constants $d_0$ and $D$ independent from $\varepsilon$. Moreover, with respect to this metric $l_{g_{\varepsilon}}(\gamma)=\tilde d_\varepsilon(\tilde x,\sigma_p\tilde x)=\frac{2\pi\varepsilon}{p}$ when $\tilde x\in\tilde\gamma$, which proves (2). Let us now prove (3).
We choose a point $x_0\in X$ such that $\inj(x,g_\varepsilon)\geq i_0>0$ (such a point exists since $g_\varepsilon=g_1$ in $Y\setminus B_y$). We define the norms 
$\parallel\gamma\parallel_\varepsilon=\tilde d_{\varepsilon}(\tilde x,\gamma (\tilde x))$ on $\Gamma=\pi_1(X)=\Z_p*G$, where $\tilde x$ is in the $x$-fiber in $\tilde X$. For all $\varepsilon$ the sets
$\Sigma_\varepsilon=\{\gamma\in\Gamma\;|\;\parallel\gamma\parallel_\varepsilon\leq 3D\}$
are generating sets for $\Z_p*G$. So:
\begin{equation}\label{Counter1}
d_1(\tilde x,\gamma\tilde x)\leq 3D\cdot\parallel\gamma\parallel_{\Sigma_1}\leq 3DS\cdot\parallel\gamma\parallel_{\Sigma_\varepsilon}\leq\frac{DS}{d_0}\cdot (3\,d_{\varepsilon}(\tilde x,\gamma\tilde x)+1)
\end{equation}
for some $S<\infty$. Actually the first and the last inequality are well known (see \cite{Gromov}, 3.22); for the second one let us define
$S^\varepsilon=\sup\{\parallel\gamma\parallel_{\Sigma_1}|\mbox{ }\gamma\in\Sigma_\varepsilon\}$ and we show that  $$S=\sup\{S^\varepsilon\mbox{ }|\mbox{ }\varepsilon\in(0,1]\}<+\infty$$ 
In fact since the sectional curvatures  (and hence the Ricci curvature) of $g_\varepsilon$ are bounded below independently from $\varepsilon$, there exists a $N(D,i_0)\in\N$ (independent from $\varepsilon$) bounding the maximum number of disjoint $g_\varepsilon$-balls of radius $i_0$ in a $g_\varepsilon$-ball of radius $3D$; it follows that $\#\Sigma_\varepsilon\leq N(D,i_0)$, for all $\varepsilon>0$; moreover  $\Sigma_\varepsilon\subseteq\Sigma_{\varepsilon'}$ for $\varepsilon'<\varepsilon$ as $d_{\varepsilon'}\leq d_{\varepsilon}$. We deduce that the sets $\Sigma_\varepsilon$ are all included in a maximal finite subset $\Sigma$, so $S<+\infty$. Then from estimate (\ref{Counter1}) we deduce readily: 
$\Entvol(X,g_\varepsilon)\leq\frac{3DS}{d_0}\Entvol(X,g_1)= H.\,\Box$
\end{exmp}

\begin{exmp}  This example shows the necessity of the upper bound for the diameter in Theorem \ref{MT2} and the 'sharpness' of the result.
Let us denote $M_1^\varepsilon=(S^1\times S^{n-1},\varepsilon^2\cdot g_0),\quad  M_2^{\varepsilon'}=(S^{1}\times S^{n-1},\frac{1}{(\varepsilon')^2}\cdot g_0)$
(where $g_0$ is the canonical product metric of $S^1\times S^{n-1}$). Let $X=M_1^\varepsilon\#M_2^{\varepsilon'}$ where the metric is constructed as follows: we cut a geodesic ball $B_1\subset M_1^{\varepsilon}$ (resp. $B_2\subset M_2^{\varepsilon'}$) of radius $\varepsilon^2$. Consider the cylinder $C=[0,1]\times S^{n-1}$; we endow $\{0\}\times S^{n-1}$ (resp. $\{1\}\times S^{n-1}$) with a Riemannian metric $h_0$ (resp. $h_1$) such that $\{0\}\times S^{n-1}$ is isometric to $(\partial B_1,\varepsilon^2\cdot g_0)$ (resp. $\{1\}\times S^{n-1}$ is isometric to $(\partial B_2,\frac{1}{(\varepsilon')^2}\cdot g_0)$). Next we define the following metric on $C$: $h_{\varepsilon,\varepsilon'}=(dr)^2+ h_r$ where $h_r$ is the metric on $\{r\}\times S^{n-1}$ defined by $h_r=(1-r)\, h_0+r\, h_1$. Finally we construct the connected sum gluing $M_1^\varepsilon$ and $M_2^{\varepsilon'}$ at the two boundary component of $C$, and we construct the metric $g_{\varepsilon,\varepsilon}$ as follows:
\begin{equation*}
g_{\varepsilon,\varepsilon'}=
\left\{
\begin{array}{c}
\varepsilon^2\cdot g_0\quad\mbox{ on }M_1^{\varepsilon}\setminus B_1;\\ 
h_{\varepsilon,\varepsilon'}\quad\mbox{ on } C;\\
\frac{1}{(\varepsilon')^2}\cdot g_0\quad\mbox{ on } M_2^{\varepsilon'}\setminus B_2;\\
\end{array}
\right.
\end{equation*}
We remark that $g_{\varepsilon,\varepsilon'}$ is not $C^\infty$ but just piecewise $C^\infty$. However it is not difficult to show that we can produce smooth metrics arbitrairely close in the sense $C^0$ to $g_{\varepsilon,\varepsilon'}$. That is why we are allowed to use the metric $g_{\varepsilon,\varepsilon'}$.\\ Let $a$, (resp. $b$) be the generator of the image of $\pi_1(M_1^\varepsilon)$ (resp. $\pi_1(M_2^{\varepsilon'})$) in the free product.
By construction $\sys\pi_1(X,g_{\varepsilon,\varepsilon'})$ is the length of the periodic geodesic freely homotopic to the geodesic loop $a$, so that:
\begin{equation}\label{sys}
\sys\pi_1(X,g_{\varepsilon,\varepsilon'})=2\pi\varepsilon
\end{equation}
On the other hand the diameter $D_{\varepsilon,\varepsilon'}$ of $(X,g_{\varepsilon,\varepsilon'})$ satisfies
\begin{equation}\label{diam}
\frac{\pi}{\varepsilon'}+1\leq D_{\varepsilon,\varepsilon'}\leq\frac{\pi}{\varepsilon'}+1+\pi\varepsilon
\end{equation}
Take $x\in\partial B_2$. Every geodesic loop based at $x$ can be written in $\pi_1(X,x)$ in the form: $\gamma=a^{p_1} b^{q_1}\cdots a^{p_m}b^{q_m}$. Such  a decomposition corresponds to a partition of the loop $\gamma$ by points $x=x_0, \,y_0,\,x_1,...,\, x_{m-1},\, y_{m-1},\,x_m=x\in \partial B_2$ such that  $a^{p_i}$ (resp. $b^{q_i}$) is the homotopy class of the loop obtained by composition of $\alpha_i$ (resp. $\beta_i$), the portion of the path $\gamma$ corresponding to $[x_{i-1},y_{i-1}]$ (resp. $[y_{i-1},x_i]$), with the minimizing geodesics joining $x$ with $x_{i-1}$ and $y_{i-1}$ (resp. with $y_{i-1}$ and $x_i$) in $\partial B_2$, whose lengths are bounded above by $\varepsilon^2\cdot C$ where $C\leq 2\pi+1$. Hence we find:
$$l(\gamma)\geq\sum_{i=1}^m(l(\alpha_i)-2C\varepsilon^2)+\sum_{i=1}^m(l(\beta_i)-2C\varepsilon^2)$$
By construction we have
$l(\alpha_i)\geq(2\pi\varepsilon)\,|p_i|+1,\, l(\beta_i)\geq\frac{2\pi}{\varepsilon'}\,|q_i|,$
so that taking $\varepsilon$ sufficiently small we get:
$l(\gamma)\geq\sum_{i=1}^m|p_i|\,(2\pi\varepsilon)+\sum_{i=1}^m|q_i|\,\frac{2\pi}{\varepsilon'}$;
hence if $\tilde x$ is in the fiber of $x$ in the Riemannian universal covering $\tilde X$ we see that
$d_{geo}(\id,\gamma)=d_{g_{\varepsilon,\varepsilon'}}(\tilde x,\gamma\tilde x)\geq d_l(\id,\gamma)$
where $d_l$ is the distance on $\Z*\Z$ corresponding to the choice of the generating system $\{a,\, b\}$ with lengths $l(a)=2\pi\varepsilon$, $l(b)=\frac{2\pi}{\varepsilon'}$. Since $X$ is compact we have:
$$\Ent(X,g_{\varepsilon,\varepsilon'})\leq\Ent(\Gamma,d_{geo})\leq\Ent(\Gamma,d_l)=h$$
where, by Lemma \ref{entsemigroup}, $h$ satisfies the equation:
\begin{equation}\label{mainex55}
(e^{2\pi h\varepsilon}-1)(e^{\frac{2\pi h}{\varepsilon'}}-1)=4
\end{equation}
\end{exmp}

5.5.1.\textit{ End of the counterexample.} If $\varepsilon=\varepsilon'$ by the estimates (\ref{sys}), (\ref{diam}) the systole and the diameter of $(X,g_{\varepsilon,\varepsilon'})$ tends respectively to $0$ and $+\infty$. On the other hand equation (\ref{mainex55}) shows that $\Ent(X,g_{\varepsilon,\varepsilon'})$ is bounded above by $\frac{1}{\pi}$. This proves the necessity of the boundedness of the diameter in Theorem \ref{MT2}.\medskip\\

5.5.2.\textit{ Optimality.} By (\ref{mainex55}) we know that $H_{\varepsilon,\varepsilon'}=\Ent(X,g_{\varepsilon,\varepsilon'})$ satisfies
$$2\pi\varepsilon\leq\frac{1}{H_{\varepsilon,\varepsilon'}}\cdot\log\left(1+\frac{4}{e^{2\cdot H_{\varepsilon,\varepsilon'}\cdot\frac{\pi}{\varepsilon'}}-1}\right)$$
Since $\frac{\pi}{\varepsilon'}\simeq\diam(X,g_{\varepsilon,\varepsilon'})$ the estimate given of Theorem \ref{MT2} is optimal.

\section*{Appendix}

This appendix is devoted to the proof of the Lemma \ref{PPFP} that we used in the proof of Theorem \ref{MT1}. We recall that a word $\gamma=\alpha_1\cdots\alpha_p$ in $\Gamma=A*B$ is said to be in the reduced form if $\forall i$ $\alpha_i\in A$ or $\alpha_i\in B$ and if $\alpha_i, \alpha_{i+1}$ do not belong to the same factor in $\Gamma$ for all $i=1,..,p-1$; in this case the length of the reduced word is $l(\gamma)=p$. We remark that this corresponds to the 'algebraic length' of $\Gamma$ only if we consider $A\sqcup B$ as the generator system of $A*B$. Notice that the reduced form is unique. We say that a word $\gamma$ is cyclically reduced if its reduced form $\gamma=\alpha_1\cdots\alpha_p$ is such that $\alpha_p\neq\alpha_1^{-1}$. \newline 

\textbf{Proof of Lemma \ref{PPFP}}  Let $\gamma=\alpha_1\cdots\alpha_p$, $\gamma'=\alpha_1'\cdots\alpha'_{p'}$ be the reduced forms for $\gamma$, $\gamma'$.
\begin{itemize}
\item If $p$ is even, then $\alpha_1$, $\alpha_p$  belong to two different factors, and the reduced form for $\gamma^r$ is: 
$$\gamma^r=(\alpha_1\cdots\alpha_p)\cdots(\alpha_1\cdots\alpha_p)$$
and the initial letter in the reduced form is not the inverse of the final one, \textit{i.e} $\gamma^r$ is cyclically reduced. It is clear that, in this case, the knowledge of $p$ and $\gamma^r$ allows us to recover the whole sequence of letters $\alpha_1,..,\alpha_p$.
\item If $p$ is odd, then $\alpha_1, \alpha_p$ belong to the same factor and the writing
$$\gamma^r=(\alpha_1\cdots\alpha_p)\cdots(\alpha_1\cdots\alpha_p)$$
can be reduced a first time by grouping together $(\alpha_p\alpha_1)$; we can not reduce further unless $\alpha_p=\alpha_1^{-1}$, and so on until $\alpha_{p-i}\neq(\alpha_{i+1})^{-1}$ this condition beeing realized for some $i\leq[\frac{p}{2}]$.  We find
$$\gamma=(\alpha_1\cdots\alpha_i)(\alpha_{i+1}\cdots\alpha_{p-i})(\alpha_{i}^{-1}\cdots\alpha_1^{-1})$$
and we get that $\gamma^r=(\alpha_1\cdots\alpha_i)(\alpha_{i+1}\cdots\alpha_{p-i})^r(\alpha_i^{-1}\cdots\alpha_1^{-1})$, so that, grouping together $\alpha_{p-i}$ and $\alpha_{i+1}$, we obtain $l(\gamma^r)\leq pr-(2i+1)(r-1)$. Moreover the initial letter in the reduced form of $\gamma^r$ is the inverse of the final one. Also in this case the knowledge of $\gamma^r$ impose the values of $\alpha_1,...,\alpha_p$.
\end{itemize}
The same arguments hold for the decomposition of $\gamma'$, \textit{i.e.} given $(\gamma')^{r'}$, $\alpha_1',..,\alpha_{p'}'$ are determined. In particular $\gamma^r=(\gamma')^{r'}$ implies  that the initial and the final letter in the reduced forms of $\gamma^r$, $(\gamma')^{r'}$ are the same. Thus the first (resp. the final) letter of the reduced word corresponding to $\gamma$ lies in the same subgroup ($A$ or $B$) of the first (resp. the final) letter of $\gamma'$; this implies that $p, p'$ are both even or both odd.
So we are led to consider the following cases:
\begin{case}[1] $p, q$ \textit{even}.
Let $w$ be the word in the alphabet $A^*\sqcup B^*$ given by the reduced form of $\gamma^r=(\gamma')^{r'}$ above; since $w$ is invariant by the shift of $p$ and $p'$ places, then it is also invariant by the shift of $d=\GCD(p,p')$ places. Therefore, setting $\tau=\alpha_1\cdots\alpha_d$ we have $\gamma=\tau^q$, $\gamma'=\tau^{q'}$ for $q=p/d$ and $q'=p'/d$.
\end{case}

\begin{case}[2] $p, q$ \textit{odd}. We know that
$$\gamma^r=(\alpha_1\cdots\alpha_i)(\alpha_{i+1}\cdots\alpha_{p-i})^r(\alpha_1...\alpha_i)^{-1}$$
$$(\gamma')^{r'}=(\alpha_1'\cdots\alpha_{i'}')(\alpha_{i'+1}'\cdots\alpha_{p'-i'}')^{r'}(\alpha_1'\cdots\alpha_{i'}')^{-1}$$
with $\gamma_1=(\alpha_{i+1}\cdots\alpha_{p-i})$ and $\gamma_1'=(\alpha_{i'+1}'\cdots\alpha_{p'-i'}')$ cyclically reduced hence, comparing the two expressions we deduce that $i=i'$ and $\alpha_k=\alpha_k'$ for $k\leq i$. Now consider $\gamma_1^r=(\alpha_{i+1}\cdots\alpha_{p-i})^r$, $(\gamma_1')^{r'}=(\alpha_{i+1}\cdots\alpha_{p'-i})^{r'}$. We have $\gamma_1^r=(\gamma_1')^{r'}$. As $\gamma_1$, $\gamma_1'$ are cyclically reduced, with $l(\gamma_1)$, $l(\gamma_1')$ odd, the only reduction that we can perform on $\gamma_1^r$ is to group together $\alpha_{p-i}$ and $\alpha_{i+1}$ (and $\alpha_{p'-i}$, $\alpha_{i'+1}$ in $(\gamma_1')^{r'}$):
$$\gamma_1^{r}=\alpha_{i+1}\cdots(\alpha_{p-i}\alpha_{i+1})\cdots(\alpha_{p-i}\alpha_{i+1})\alpha_{i+2}\cdots\alpha_{p-i}$$
$$(\gamma_1')^{r}=\alpha_{i+1}'\cdots(\alpha_{p'-i}'\alpha_{i+1}')\cdots(\alpha_{p'-i}'\alpha_{i+1}')\alpha_{i+2}'\cdots\alpha_{p'-i}'$$
this implies $\alpha_{i+1}=\alpha_{i+1}'$ and so setting
$$\tilde\gamma_1=\alpha_{i+2}\cdots\alpha_{p-i-1}(\alpha_{p-i}\alpha_{i+1}), \quad\tilde\gamma_1'=\alpha_{i+1}'\cdots\alpha_{p'-i-1}'(\alpha_{p'-i}'\alpha_{i+1}')$$
we have $\tilde\gamma_1^r=(\alpha_{i+1})^{-1}\,\gamma_1^r\,\alpha_{i+1}=(\alpha_{i+1})^{-1}\,(\gamma_1')^{r'}\,\alpha_{i+1}=(\tilde\gamma_1')^{r'}$ and we are reduced to the case where $l(\tilde\gamma_1)$ and $l(\tilde\gamma_1')$ are even, which we treat as before. Therefore we can find a $\tilde\tau$ and integers $q, q'$ such that $\tilde\gamma_1=\tilde\tau^q$, $\tilde\gamma_1'=(\tilde\tau')^{q'}$. Setting $\tau=\alpha\tilde\tau\alpha^{-1}$ for $\alpha=\alpha_1\cdots\alpha_{i+1}$ we finally have: $\gamma=\tau^q$, $\gamma'=\tau^{q'}$. $\Box$

\end{case}

\end{document}